\documentclass[12pt,A4]{article}
\usepackage{amssymb, amsmath, enumerate}
\begin{document}
\newtheorem{thm}{Theorem}
\newtheorem{lm}[thm]{Lemma}
\newtheorem{prop}[thm]{Proposition}
\newtheorem{cor}[thm]{Corollary}
\newtheorem{defi}[thm]{Definition}
\newtheorem{eg}[thm]{Example}

\numberwithin{equation}{section}

\newcommand{\proof}{{\bf Proof\ }}
\newcommand{\eproof}{$\blacksquare$\bigskip}
\newcommand{\CC}{\mathcal C}
\newcommand{\CF}{\mathcal F}
\newcommand{\CR}{\mathcal R}
\newcommand{\C}{\Bbb C}
\newcommand{\R}{\Bbb R}
\newcommand{\Z}{\Bbb Z}
\newcommand{\N}{\Bbb N}

\newcommand{\prob}{\mathop{\mathrm{Pr}}}
\newcommand{\eps}{\varepsilon}
\newcommand{\raction}{\triangleleft}
\newcommand{\Del}{\Delta}
\newcommand{\del}{\delta}
\newcommand{\Lan}{\Lambda}
\newcommand{\Gam}{\Gamma}
\newcommand{\gam}{\gamma}
\newcommand{\ten}{\bigotimes}
\newcommand{\sten}{\otimes}
\newcommand{\Om}{\Omega}
\newcommand{\sfi}{\varphi}
\newcommand{\om}{\omega}
\newcommand{\lan}{\lambda}
\newcommand{\sig}{\sigma}
\newcommand{\al}{\alpha}
\newcommand{\be}{\beta}
\newcommand{\all}{\forall}
\newcommand{\flsh}{\longrightarrow}
\newcommand{\ld}{\ldots}
\newcommand{\cd}{\cdots}
\newcommand{\bu}{\bullet}
\newcommand{\p}{\partial}
\newcommand{\eqn}[2]{\begin{equation}#2\label{#1}\end{equation}}

\textheight 23.6cm \textwidth 16cm \topmargin -.2in \headheight 0in
\headsep 0in \oddsidemargin 0in \evensidemargin 0in \topskip 28pt

\title{On the Series\\ $\sum_{n=0}^{\infty}\frac{(-1)^n}{n!}t^nf^{(n)}(t)$}

\author{S.E.Akrami\footnote{Supported financially by the grant 83810319 from Institute for Research in Fundamental
Sciences, P.O.Box:19395-5746, Tehran, Iran.
E-mail:akramisa@ipm.ir}\\
Departments of Mathematics and Physics,\\ Institute for Research in
Fundamental Sciences,\\ P.O.Box:19395-5746, Tehran, Iran} \maketitle
\begin{abstract}
We study the series
$\sum_{n=0}^{\infty}\frac{(-1)^n}{n!}t^nf^{(n)}(t)$. We show that
for analytic functions this series is uniformly and absolutely
convergent to the constant $f(0)$. We show that there are nowhere
analytic functions for them the series is divergent for all $t$ and
also there are nowhere analytic functions for them the series is
convergent to $f(0)$ at least for $t$ in a dense subset of $\R$.
\end{abstract}
\textbf{Acknowledgement} I thank Allah Rabbel-Alamin and Imam Zaman.
\section{Introduction}

This paper is about the convergence of the series
\eqn{hat0}{\hat{f}(t):=\sum_{n=0}^{\infty}(-1)^n\frac{f^{(n)}(t)}{n!}t^n.}

First notice that if we differentiate the series term by term we get
\begin{eqnarray}\frac{d}{dt}\hat{f}(t)&=&
\sum_{n=0}^{\infty}(-1)^n\frac{f^{(n+1)}(t)}{n!}t^n+
\sum_{n=1}^{\infty}(-1)^n\frac{f^{(n)}(t)}{(n-1)!}t^{n-1}\nonumber\\&=&
\sum_{n=0}^{\infty}(-1)^n\frac{f^{(n+1)}(t)}{n!}t^n-\sum_{n=0}^{\infty}(-1)^n\frac{f^{(n+1)}(t)}{n!}t^n
\nonumber\\&=&0\nonumber
\end{eqnarray}
Thus $\hat{f}(t)$ should be constant. But the point is that we are
not allowed to differentiate term by term from a series even the
series is uniformly convergent. We shall show that for analytic
functions around origin this series is convergent to the constant
$f(0)$. The surprising point here is not the proof which is very
easy, but it is very strange that why this fact has been forgotten
in textbooks of mathematical analysis! At least it might be
mentioned in them as an exercise. However I could not find any trace
of this strange series in the mathematics literature, I could see it
in quantum mechanics textbooks, see for example \cite{BL}. There, it
is named as the translation operator. Physicists define an operator
as \eqn{T}{(T_af)(t):=f(t+a).} Then they claim that this operator is
equal to the operator
\eqn{e}{e^{a\frac{d}{dt}}=\sum_{n=0}^\infty\frac{a^n}{n!}\frac{d^n}{dt^n}.}But
in fact to prove this they implicitly or sometimes explicitly assume
that $f$ is analytic. In fact one can easily prove, as we shall
prove in next section, that for analytic functions these two
operators coincide. Notice that if we set $a=-t$ we get our strange
series.

Some textbooks,\cite{BL}, try to prove the equality of the operators
(\ref{T}) and (\ref{e}) as follows. They first expand $f$ as a
Fourier integral \eqn{}{f(t)=\int_{-\infty}^{\infty}a(x)e^{ixt}dx,}
and then apply the operator (\ref{e}) to this integral. They
exchange the order of the derivative $\frac{d^n}{dt^n}$ with the
integral.
\begin{eqnarray}e^{-t\frac{d}{dt}}f&=&\sum_{n=0}^\infty\frac{(-1)^n}{n!}t^nf^{(n)}(t)\nonumber\\&=&
\sum_{n=0}^\infty\int_{-\infty}^{\infty}\frac{(-1)^n}{n!}(ix)^nt^na(x)e^{ixt}dx
\nonumber\\&=&
\int_{-\infty}^{\infty}\sum_{n=0}^\infty\frac{(-ixt)^n}{n!}a(x)e^{ixt}dx\nonumber\\&=&
\int_{-\infty}^{\infty}e^{-ixt}a(x)e^{ixt}dx\nonumber\\&=&\int_{-\infty}^{\infty}a(x)dx\nonumber\\&=&f(0)\nonumber
\end{eqnarray}
But in the third line of the above proof we have exchanged the order
of integral and summation, but this is not allowed unless the series
$\sum_{n=0}^\infty\frac{(-ixt)^n}{n!}a(x)e^{ixt}=a(x)$ is uniformly
convergent. But this is not always true and depends to the
coefficient $a(x)$. For example if $a(x)$ has compact support then
the series is uniformly convergent.

\section{Convergence Tests}
\begin{prop}
A necessary condition for the point-wise (uniformly) convergence of
the series $\hat{f}$ is that
\eqn{}{\lim_{n\rightarrow\infty}\frac{1}{n!}t^nf^{(n)}(t)=0,}
point-wise (uniformly).
\end{prop}
\proof A necessary condition for the convergence of a series $\sum_n
a_n$ is $\lim_{n\rightarrow\infty}a_n=0$.\eproof
\begin{thm}\label{main condition}
Let for some $k>0$ the limit
\eqn{main}{\lim_{n\rightarrow\infty}\frac{(-1)^n}{n!}t^nf^{(n+k)}(t)}
exists uniformly around origin. If $k>1$ then the series $\hat{f}$
is uniformly convergent to the constant $f(0)$ and if $k=1$ then the
series $\hat{f}$ is uniformly convergent and differentiable. If when
$k=1$ we denote the limit (\ref{main}) by $\check{f}(t)$, then we
have $\frac{d}{dt}\hat{f}(t)=\check{f}(t)$. In particular case when
$\check{f}(t)=0$ for all $t$, i.e.
\eqn{main2}{\lim_{n\rightarrow\infty}\frac{(-1)^n}{n!}t^nf^{(n+1)}(t)\rightarrow0}
uniformly, we get $\hat{f}(t)=f(0).$
\end{thm}
\proof  Let $S_N$ be the partial sum of the series obtained by the
term by term differentiation of the series $\hat{f}$. Thus
\begin{eqnarray}S_N&=&\sum_{n=0}^N\frac{(-1)^n}{n!}t^nf^{(n+1)}(t)+
\sum_{n=1}^{N}\frac{(-1)^n}{(n-1)!}t^{n-1}f^{(n)}(t)\nonumber\\&=&\sum_{n=0}^N\frac{(-1)^n}{n!}t^nf^{(n+1)}(t)-
\sum_{n=0}^{N-1}\frac{(-1)^n}{n!}t^nf^{(n+1)}(t)
\nonumber\\&=&\frac{(-1)^N}{N!}t^Nf^{(N+1)}(t).\nonumber
\end{eqnarray}
Now let us define
\eqn{}{f_{k,n}(t):=\frac{(-1)^n}{n!}t^nf^{(n+k)}(t).} We have
$f_{k-1,n}'(t)=f_{k,n}(t)-f_{k,n-1}(t)$. Thus if for some $k>1$,
$\lim_{n\rightarrow\infty}f_{k,n}(t)$ exists uniformly then
$\lim_{n\rightarrow\infty}f_{k-1,n}'(t)=0$ uniformly. Thus
$\lim_{n\rightarrow\infty}f_{k-1,n}(t)$ exists uniformly and is
constant and since we have $f_{k,n}(0)=0$ for all $k$ we deduce that
$\lim_{n\rightarrow\infty}f_{k-1,n}(t)=0$ uniformly. By repeating
this argument we deduce that $\lim_{n\rightarrow\infty}f_{1,n}(t)=0$
uniformly. But we have $f_{1,N}(t)=S_N$. Thus the series obtained by
the term by term differentiation of the series $\hat{f}$ is
uniformly convergent and since the series $\hat{f}$ is convergent at
$t=0$, then by a well known theorem in mathematical analysis (see
\cite{Apo} for example) the series $\hat{f}$ is also uniformly
convergent and differentiable and
$\frac{d}{dt}\hat{f}(t)=\lim_{N\rightarrow\infty}S_N=0$. Thus
$\hat{f}(t)=f(0)$.

Now if $k=1$ then
$\lim_{N\rightarrow\infty}S_N=\lim_{N\rightarrow\infty}f_{1,N}(t)$
exists uniformly but not necessarily vanishes. Thus the series
obtained by the term by term differentiation of the series $\hat{f}$
is uniformly convergent and therefore the series $\hat{f}$ is also
uniformly convergent and differentiable and
$\frac{d}{dt}\hat{f}(t)=\check{f}(t)$.\eproof
\begin{thm}\label{power}
If there exists constants $C$ and $M$ such that
\eqn{}{|f^{(n)}(t)|<CM^n} around origin. Then the series $\hat{f}$
is uniformly convergent to the constant $f(0)$.
\end{thm}
\proof Let $t\in(-a,a)$ for some $a$. Then we have
$|\frac{(-1)}{n!}t^nf^{(n+1)}(t)|\le CM\frac{(Ma)^n}{n!}$. But it is
easy to show that the right hand side of this inequality goes to
zero. Thus we can use the theorem \ref{main condition}. \eproof

\begin{thm}
If there exists a constant $M$ such that
$|\frac{f^{(n+1)}(t)}{f^{(n)}(t)}|<M$  around origin. Then the
series $\hat{f}$ is uniformly convergent to the constant $f(0)$.
\end{thm}
\proof Let $t\in[-a,a]$ for some $a$. We have
$|f^{(n+1)}(t)|<M|f^{(n)}(t)|$. Thus $|f^{(n)}(t)|<M^n|f(t)|$. Now
let $C$ be the maximum of $f$ at $[-a,a]$. Thus $|f^{(n)}(t)|<CM^n.$
Now use the theorem \ref{power}. \eproof
\begin{thm}
If there exists a constant $M$ such that $|\sum_{n=0}^N
(-1)^nf^{(n)}(t)|<M$ for all $N$  uniformly around origin. Then the
series $\hat{f}$ is uniformly convergent but not necessarily to a
constant.
\end{thm}
\proof We use the Dirichlet' test (see \cite{Apo} for example) which
states that if there exists a constant $M$ such that $|\sum_{n=0}^N
f_n(t)|<M$ for all $N$  uniformly and $g_n(t)$ is a decreasing
sequence converging to zero uniformly, then the series $\sum
f_n(t)g_n(t)$ is uniformly convergent. In our case we set
$f_n(t):=(-1)^nf^{(n)}(t)$ and $g_n(t):=t^n/n!$.\eproof
\begin{thm}
Let
\eqn{}{\al:=\sup_{t}\limsup_{n}\sqrt[n]{\frac{|f^{(n)}(t)|}{n!}},~~~
\be:=\inf_{t}\limsup_{n}\sqrt[n]{\frac{|f^{(n)}(t)|}{n!}},} and
$R:=\frac{1}{\al}, S:=\frac{1}{\be}.$ Then the series $\hat{f}(t)$
for $|t|<R$ is absolutely convergent but not necessarily to $f(0)$
and for $|t|>S$ is divergent.
\end{thm}
\proof Let $a_n(t):=\frac{(-1)}{n!}f^{(n)}(t)$. Thus
$\hat{f}(t)=\sum_{n=0}^{\infty}a_n(t)t^n$. Now for $|t|<R$ we have
$\limsup_{n}\sqrt[n]{|a_n(t)t^n|}=|t|\limsup_{n}\sqrt[n]{\frac{|f^{(n)}(t)|}{n!}}\le\al|t|<1$.
Thus the series $\sum_{n=0}^{\infty}a_n(t)t^n=\hat{f}(t)$ is
absolutely convergent. The other part is similar.\eproof
\begin{prop}\label{linear}
(i)If the series $\hat{f}$ and $\hat{g}$ are either point-wise or
uniformly and either absolutely or conditional convergent then the
series $\widehat{af+bg}$ are so and
$\widehat{af+bg}=a\hat{f}+b\hat{g}$

(ii) If at least one of series $\hat{f}$ and $\hat{g}$ are
absolutely (either point-wise or uniformly) convergent then the
series $\widehat{fg}$ is also (either point-wise or uniformly)
convergent and $\widehat{fg}=\hat{f}\hat{g}$.
\end{prop}
\proof (i) is obvious. (ii) Using the fact that
$(fg)^{(n)}=\sum_{i=0}^n\frac{n!}{i!(n-i)!}f^{(i)}g^{(n-i)}$ we have
\begin{eqnarray}\hat{f}(t)\hat{g}(t)&=&
\sum_{n=0}^{\infty}(-1)^n\frac{f^{(n)}(t)}{n!}t^n\sum_{n=0}^{\infty}(-1)^n\frac{g^{(n)}(t)}{n!}t^n\nonumber\\&=&
\sum_{n=0}^{\infty}\sum_{i+j=n}(-1)^i\frac{f^{(i)}(t)}{i!}(-1)^j\frac{g^{(j)}(t)}{j!}t^it^j\nonumber\\&=&
\sum_{n=0}^{\infty}\sum_{i=0}^n\frac{n!}{i!(n-i)!}f^{(i)}(t)g^{(n-i)}(t)\frac{(-1)^n}{n!}t^n\nonumber\\&=&
\sum_{n=0}^{\infty}(-1)^n\frac{(fg)^{(n)}(t)}{n!}t^n\nonumber\\&=&\widehat{fg}(t)\nonumber
\end{eqnarray}
Where we used the Mertens' theorem about the convergence of the
Cauchy's product of two series \cite{Apo}.\eproof
\begin{prop}\label{scale}
If $g(t):=f(at)$ for some constant $a$. Then
$\hat{g}(t)=\hat{f}(at)$.
\end{prop}
\proof
$\hat{g}(t):=\sum_{n=0}^{\infty}(-1)^n\frac{g^{(n)}(t)}{n!}t^n=
\sum_{n=0}^{\infty}(-1)^n\frac{f^{(n)}(at)}{n!}a^nt^n=\hat{f}(at)$.
\eproof

\begin{thm}
If the series $\hat{f}$ is uniformly convergent (not necessarily to
a constant) then for the anti-derivative of $f$, i.e.
$F(x)=\int_0^xf(t)dt$ the series $\hat{F}$ is uniformly convergent
to the constant $F(0)=0$.
\end{thm}
\proof First, notice that if we set
$f_n(x):=\int_0^xt^nf^{(n)}(t)dt$ then by the integration by parts
we get $f_n(x):=x^nf^{(n-1)}(x)-nf_{n-1}(x)$. Thus
$x^nF^{(n)}(x)=x^nf^{(n-1)}(x)=f_n(x)+nf_{n-1}(x)$. Now since
$\hat{f}$ is a uniformly convergent series of integrable functions,
we can integrate term by term, \cite{Apo}. Thus
\begin{eqnarray}\int_{0}^x\hat{f}(t)dt&=&
\sum_{n=0}^{\infty}\frac{(-1)^n}{n!}\int_0^xt^nf^{(n)}(t)dt\nonumber\\&=&
\sum_{n=0}^{\infty}\frac{(-1)^n}{n!}f_n(x)\nonumber
\end{eqnarray}
Next we get
\begin{eqnarray}\hat{F}(x)&=&F(x)+\sum_{n=1}^{\infty}\frac{(-1)^n}{n!}x^nF^{(n)}(x)\nonumber\\&=&F(x)+
\sum_{n=1}^{\infty}\frac{(-1)^n}{n!}f_n(x)-\sum_{n=1}^{\infty}\frac{(-1)^{n-1}}{(n-1)!}f_{n-1}(x)
\nonumber\\&=&F(x)+\sum_{n=1}^{\infty}\frac{(-1)^n}{n!}f_n(x)-\sum_{n=0}^{\infty}\frac{(-1)^{n}}{n!}f_{n}(x)
\nonumber\\&=&F(x)-\int_0^xf(t)dt\nonumber\\&=&0\nonumber
\end{eqnarray}\eproof
\begin{thm}\label{analytic}
If a function $f$ can be expanded around origin via power series
$f(t)=\sum_{n=0}^{\infty}\frac{f^{(n)}(0)}{n!}t^n,|t|<R$, then
$\hat{f}$ is absolutely and uniformly convergent to constant $f(0)$
in the interval $|t|<R/2$. That is
\eqn{hat0}{\sum_{n=0}^{\infty}(-1)^n\frac{f^{(n)}(t)}{n!}t^n=f(0),~~~|t|<R/2,}
for analytic functions. If $R=\infty$ then $\hat{f}(t)=f(0)$for all
$t\in\R$.
\end{thm}
\proof It is known that (see \cite{Apo} for example) for any $x$
satisfying $|x|<R$, we can expand $f$ around $x$ as
$f(t)=\sum_{n=0}^{\infty}\frac{f^{(n)}(x)}{n!}(t-x)^n$. This holds
for all $t$ satisfying $|t-x|<R-|x|$. Now if $|x|<R/2$ then
$|0-x|<R-|x|$. Thus we can put $x=0$ in the last series to get
$f(0)=\sum_{n=0}^{\infty}\frac{f^{(n)}(x)}{n!}(-x)^n=\hat{f}(x)$.
\eproof

The following example shows that we can not make larger the domain
of validity of $\hat{f}(x)=f(0)$ from the domain $|x|<R/2$.
\begin{eg}\label{rational}
For the function $f(t)=\frac{1}{1+t}$ the series $\hat{f}$ is
absolutely and uniformly convergent to the constant $f(0)$ for
$\frac{-1}{2}<t$ and otherwise is divergent.
\end{eg}
\proof Using induction on $n$ we have
$f^{(n)}(t)=(-1)^{n}n!\frac{1}{(1+t)^{n+1}}$. Thus
\begin{eqnarray}\hat{f}(t)&=&
\sum_{n=0}^{\infty}\frac{t^n}{(1+t)^{n+1}}\nonumber\\
&=&\frac{1}{1+t}\sum_{n=0}^{\infty}(\frac{t}{1+t})^n\nonumber\\&=&
\frac{1}{1+t}\frac{1}{1-\frac{t}{1+t}}\nonumber\\&=&1\nonumber\\&=&f(0).\nonumber
\end{eqnarray} The convergence holds  when $|\frac{t}{1+t}|<1.$ i.e when $\frac{-1}{2}<t$.\eproof
\begin{thm}\label{analytic2}
Let $f$ be a smooth function. Let $t\ne0$ belongs to its domain.
Consider the Taylor series around $t$
\eqn{Tf}{(T_tf)(x):=\sum_{n=0}^{\infty}\frac{f^{(n)}(t)}{n!}(x-t)^n,}
whose radius of convergence is denoted by $R_t$.

(i) If $0$ does not belong to the convergence interval of
(\ref{Tf}), i.e. if $R_t<|t|$, in particular if $R_t=0$, then
$\hat{f}(t)$ diverges.

(ii) If $0$  belongs to the convergence interval of (\ref{Tf}), i.e.
if $R_t>|t|$, then $\hat{f}(t)$ converges to $(T_tf)(0)$. Moreover
if  $f$ is analytic at $t$, i.e. if $T_tf=f$ then $\hat{f}(t)$
converges to $f(0)$.

(iii) If $0$  belongs to the boundary of the convergence interval of
(\ref{Tf}), i.e. if $R_t=|t|$, then in general one can not say
anything about the convergence of $\hat{f}(t)$. But if we know
priory that $(T_tf)(x)$ converges at $x=0$ or equivalently
$\hat{f}(t)$ is convergent then the limit $\lim_{x\rightarrow
0}(T_tf)(x)$ exists and is equal to $\hat{f}(t)$. In particular if
$f$ is analytic at $t$ and moreover $0$ belongs to the domain of $f$
then $\hat{f}(t)=f(0)$.
\end{thm}
\proof Since indeed $\hat{f}(t)=(T_tf)(0)$ the cases $(i)$ and
$(ii)$ are clear. Case (iii) is the Able's theorem \cite{Apo}.
\eproof

\textbf{Remark} 1)The example $f(t)=\frac{1}{1+t}$ which was
investigated above using theorem (\ref{analytic}) can also be
studied by the theorem (\ref{analytic2}) with the same results.

2) This theorem gives us another proof for the theorem
(\ref{analytic}).

Now we study the series $\hat{f}(t)$ for some examples of
nonanalytic functions. The first example is about a function which
is analytic everywhere except at $t=0$
\begin{eg}
For the function \eqn{}{f(t):=e^{\frac{-1}{t^2}},~t\ne0,~f(0)=0.} If
$\hat{f}(t)$ is convergent then we must have $\hat{f}(t)=f(0)$. But
unfortunately until the time of writing of this article we do not
know anything about the convergence of $\hat{f}(t)$ for $t\ne0$.
\end{eg}
\proof If $f$ is analytic at some $t\ne0$ in which $R_t<|t|$ then by
a well known theorem of mathematical analysis \cite{Apo}, $f$ must
be analytic at  $0$. But it is a well known fact that $f$ is not
analytic at $0$. Thus we must have $R_t\ge|t|$. In fact this
function satisfies the case (iii). That is $f$ is analytic at any
$t\ne0$ and we have $R_t=|t|$. Because this function is the
composition of a everywhere analytic function $g(t)=e^t$ and the
function $h(t)=\frac{-1}{t^2}$ which is analytic everywhere except
at $t=0$ and we know that for such a composite $f=g\circ h$ function
in which $g$ is analytic everywhere, the radius of convergence of
$f$ at $t$ is equal to the radius of convergence of $h$ at $t$. And
clearly the radius of convergence of $h$ at $t\ne0$ is $R_t=|t|$.
\eproof

\begin{eg}\label{u}
Let $u$ be a smooth periodic positive function with period $\ell$,
analytic at each point except at points $m\ell$ where $m\in\Z$, for
each $n$ there exists a constant $M_n$ such that $|u^{(n)}(x)|<M_n$
for all $x$ and $u^{(n)}(0)=0$ for all $n$. Let $a_n$ be a sequence
of positive real numbers such that the power series
$\sum_{n=0}^\infty a_nx^n$ converges all over $\R$. Then the
following function \eqn{series}{f(x)=\sum_{n=0}^\infty a_nu(2^nx),}
is smooth but nowhere analytic. That is at each point $a$ the Taylor
series at $a$ either diverges or does not converge to $f(x)$.

Moreover $\hat{f}(t)$ for all
$t\in\{\frac{(2m+1)\ell}{2^n}~|~m\in\Z,n\in\N\}$ is convergent to
$f(0)=0$. In other points we do not know anything on the convergence
of $\hat{f}(t)$ and its values.
\end{eg}
\proof Since $|a_nu(2^nx)|<a_nM_0$ and since the series $\sum_n a_n$
converges, the series (\ref{series}) is uniformly convergent and
thus $f$ is well defined. Moreover we have
$|\frac{d^k}{dx^k}(a_nu(2^nx))|\le|2^{kn}a_nu^{(k)}(2^nx)|<2^{kn}a_nM_k$
and since the series $\sum_n a_n (2^k)^n$  is convergent by
hypothesis, we conclude that the series
$\sum_n\frac{d^k}{dx^k}(a_nu(2^nx))$ is uniformly convergent and
thus $f$ is smooth.

Next suppose that $f$ is analytic at some point. Then since
analyticity at a point implies analyticity at some neighborhood of
that point and since the set of numbers of the form
$\frac{m\ell}{2^n}$ where $m$ is an odd integer and $n$ is a natural
number, is dense in $\R$, we may assume that $f$ is analytic at
$a=\frac{m_0\ell}{2^{n_0}}$ for some $m_0$ and $n_0$. Now let
$g(x):=\sum_{n=0}^{n_0} a_nu(2^nx),~h(x):=\sum_{n=n_0+1}^\infty
a_nu(2^nx)$. Thus $h=f-g$ and since $g$ is analytic at $a$ then
analyticity of $f$ at $a$ implies analyticity of $h$ at $a$. But we
have $h^{(k)}(a)=\sum_{n=n_0+1}^\infty
2^{kn}a_nu^{(k)}(2^n\frac{m_0\ell}{2^{n_0}})=\sum_{n=n_0+1}^\infty
2^{kn}a_nu^{(k)}(2^{n-n_0}m_0\ell)=0$, since periodicity of $u$
implies that $u^{(k)}(2^{n-n_0}m_0\ell)=0$. Thus $h$ should vanishes
around $a$ but this is a contradiction by hypothesis of positiveness
of the sequence $a_n$ and the function $u$.

We showed that at each point of the form $a=\frac{m_0\ell}{2^{n_0}}$
for some odd integer $m_0$ and natural number $n_0$, the Taylor
series at $a$ converges to $g(x)$ for all $x\in\R$. Thus by the
theorem (\ref{analytic2}) part (ii),
$\hat{f}(\frac{m_0\ell}{2^{n_0}})$ is well defined and is equal to
$g(0)$. But we have $g(0)=f(0)=0.$\eproof

The following example is from \cite{KK}.
\begin{eg}
According to the example \ref{u}, let $a_n:=\frac{1}{n!}$ and
$u(x):=\be(x-[x]).$ Where $\be(x)$ is any smooth function on $[0,1]$
which is positive and analytic on $(0,1)$ and
$\be^{(n)}(0)=\be^{(n)}(1)=0$ for all $n$, and for each $n$ there
exists a constant $M_n$ such that $|\be^{(n)}(x)|<M_n$ for all $x$.
For instance let $\be(x):=\al(x)\al(1-x)$ where  $\al(x)$ is any
smooth function which is positive  and analytic on $(0,1)$ and
$\al^{(n)}(0)=0$ for all $n$ and for each $n$ there exists a
constant $M_n$ such that $|\al^{(n)}(x)|<M_n$ for all $x$. For
instance let $\al(x):=e^{\frac{-1}{x^s}},~\al(0)=0$ where $s=1$ or
$s=2$. Thus the functions \eqn{}{f(x)=\sum_{n=0}^\infty
\frac{1}{n!}\be(2^nx-[2^nx]),} are smooth but nowhere analytic. Thus
$\hat{f}(t)$ for all $t\in\{\frac{(2m+1)}{2^n}~|~m\in\Z,n\in\N\}$ is
convergent to $f(0)=0$. In other points we do not know anything on
the convergence of $\hat{f}(t)$ and its values.
\end{eg}
\proof We show $u$ is smooth at every point $a$. If $a$ is not an
integer then $x-[x]$ being  equal to the smooth function $x-[a]$,
around $a$,  is smooth and therefore since $\be$ is smooth, the
composition $u(x)=\be(x-[a])$ is smooth at $a$ and
$u^{(k)}(a)=\be^{(k)}(a-[a])$. Next let $a=n$ be an integer. We have
$u'(n+)=\lim_{x\rightarrow
n^+}\frac{u(x)-u(n)}{x-n}=\lim_{x\rightarrow
n^+}\frac{\be(x-n)-\be(0)}{x-n}=\be'(0)=0,$ and
$u'(n-)=\lim_{x\rightarrow
n^-}\frac{u(x)-u(n)}{x-n}=\lim_{x\rightarrow
n^-}\frac{\be(x-n+1)-\be(1)}{x-n}=\be'(1)=0$. Thus $u'(a)=0$. Now
let $u^{(k)}(n)$ exists and is equal to zero. We have
$u^{(k+1)}(n+)=\lim_{x\rightarrow
n^+}\frac{u^{(k)}(x)-u^{(k)}(n)}{x-n}=\lim_{x\rightarrow
n^+}\frac{\be^{(k)}(x-n)}{x-n}=\lim_{x\rightarrow
n^+}\frac{\be^{(k)}(x-n)-\be^{(k)}(0)}{x-n}=\be^{(k+1)}(0)=0$ and
similarly $u^{(k+1)}(n-)=\lim_{x\rightarrow
n^-}\frac{u^{(k)}(x)-u^{(k)}(n)}{x-n}=\lim_{x\rightarrow
n^-}\frac{\be^{(k)}(x-n+1)}{x-n}=\lim_{x\rightarrow
n^-}\frac{\be^{(k)}(x-n+1)-\be^{(k)}(1)}{x-n}=\be^{(k+1)}(1)=0$.
Thus $u^{(k+1)}(n)=0$. Therefore $u$ is smooth.

Clearly $u$ is periodic with period $1$ and analytic at everywhere
except at integers. Also we have
$|u^{(k)}(x)|=|\be^{(k)}(x-[x])|<M_n$.

When $\be(x):=\al(x)\al(1-x)$, by the Leibnitz rule
$$\be^{(k)}(x)=\sum_{k=0}^n
(-1)^{n-k}\frac{n!}{k!(n-k)!}\al^{(k)}(x)\al^{(n-k)}(1-x)$$ and the
fact that $\al^{(k)}(0)=0$ we conclude that
$\be^{(k)}(0)=\be^{(k)}(1)=0.$ Also $|\be^{(k)}(x)|\le\sum_{k=0}^n
\frac{n!}{k!(n-k)!}M_kM_{n-k}$. Thus derivatives of $\be$ are
bounded.

One can easily check that the functions
$\al(x):=e^{\frac{-1}{x^s}},~\al(0)=0$ where $s=1$ or $s=2$ satisfy
all the required conditions. \eproof

The following example is from \cite{BC}.
\begin{eg}
According to the example \ref{u}, let $u(x):=\al(\sin x)$ where
$\al(x):=e^{\frac{-1}{x^2}},x\ne0,~\al(0)=0$ and let
$a_n:=2^{-2^n}$. Thus the functions \eqn{}{f(x)=\sum_{n=0}^\infty
2^{-2^n}e^{-\csc^2(2^nx)},} is smooth but nowhere analytic. Thus
$\hat{f}(t)$ for all
$t\in\{\frac{2\pi(2m+1)}{2^n}~|~m\in\Z,n\in\N\}$ is convergent to
$f(0)=0$. In other points we do not know anything on the convergence
of $\hat{f}(t)$ and its values.
\end{eg}
\proof Since $u^{(n)}(x)$ is a linear combination of $\al^{(k)}(\sin
x ),1\le k\le n$ with coefficients being linear combination of sine
and cosine functions with constant coefficients and since
derivatives of $\al$ are bounded and sine and cosine are bounded
functions  we conclude that for each $n$ there exists a constant
$M_n$ such that $|u^{(n)}(x)|<M_n$ for all $x$ and $u^{(n)}(0)=0$
for all $n$. Verifying other parts of conditions of the example
\ref{u} is easy. \eproof

Another evidence for strangeness of the series $\hat{f}(t)$ comes
from the following argument. Suppose $f(t)=\sum_{m=-\infty}^\infty
c_me^{i\om _mt}.$  Then  $f^{(n)}(t)=\sum(i\om_m)^n c_me^{i\om_m
t}.$ Thus
\begin{eqnarray}\hat{f}(t)&=&\sum_n\sum_m\frac{(-1)^n(i\om_m)^{n}}{n!}c_me^{i\om_m
t}t^n\nonumber\\&=&
\sum_m\sum_n(\frac{(-1)^n(i\om_m)^{n}}{n!}t^n)c_me^{i\om_m
t}\nonumber\\&=& \sum_me^{-i\om_m t}c_me^{i\om_m
t}\nonumber\\&=&\sum_{m=-\infty}^{\infty}c_m\nonumber\\&=&f(0).\nonumber
\end{eqnarray}
But the above argument is again analytically ill, since we have
exchanged the order of two infinite sums which from the theorems of
mathematical analysis we are not allowed in general to do so. In
order to be able to do so there exists a general theorem on double
series which states that if for a double infinite series
$\sum_{m,n}a_{mn}$ for each $n$ the series $\sum_m |a_{mn}|$ is
convergent which we show its sum by $b_n$ and the series $\sum_n
b_n$ is also convergent then we are allowed to exchange the order of
summation. That is we have $\sum_n\sum_ma_{mn}=\sum_m\sum_na_{mn}$.
Now let us check if this criterion can be applied to the double
series whose entries are $a_{mn}:=\frac{(-im)^n}{n!}c_me^{im t}t^n$.
For simplicity we have assumed that $\om_m=1$ for all $m$. We have
$\sum_m |a_{mn}|=2\frac{t^n}{n!}\sum_{m=1}^\infty m^n|c_m|$. But in
Fourier analysis it is well known that the series
$f^{(n)}(t)=\sum(im)^n c_me^{im t}$ is absolutely convergent. That
is the series $\sum_{m=1}^\infty m^n|c_m|$ is convergent which we
show its sum by $\al_n$. Thus we have $\sum_m
|a_{mn}|=2\frac{\al_nt^n}{n!}.$ Thus we should verify the
convergence of  the series $\sum_{n=0}^\infty\frac{\al_nt^n}{n!}$.
But this is a power series whose convergence radius is given by
$R=\al^{-1}$ where
$\al:=\limsup\sqrt[n]{\frac{\al_n}{n!}}=\limsup\sqrt[n]{\frac{\al_n}{n!}}$.
The point is that we are not sure if $R\ne0$? See the following
examples, \cite{R}.
\begin{eg}
The function
\eqn{}{f(t):=\sum_{m=0}^{\infty}\frac{1}{m!}e^{i2^mt}}is smooth
nowhere analytic, in the sense that convergence radius of the
Taylor's series of $f$ at each point is zero and therefore
$\hat{f}(t)$ diverges for all $t\ne0$.
\end{eg}
\proof Since for all $n$ we have
$\sum_{m=0}^{N}|\frac{(2^mi)^n}{m!}e^{i2^mt}|\le\sum_{m=0}^{N}\frac{(2^mi)^n}{m!}$
and since the later series converges to $i^ne^{2^n}$, we conclude
that $f$ and its derivatives $f^{(n)}$ are well defined \cite{Apo}.
Next we have
$f^{(n)}(0)=\sum_{m=0}^{\infty}\frac{(2^mi)^n}{m!}=i^ne^{2^n}$. Thus
$\sum_{n=0}^{\infty}\frac{f^{(n)}(0)}{n!}t^n=\sum_{n=0}^{\infty}\frac{i^ne^{2^n}}{n!}t^n$.
The convergence radius is obtained by the ratio test as follows.
$|\frac{\frac{i^{n+1}e^{2^{n+1}}}{(n+1)!}}{\frac{i^ne^{2^n}}{n!}}|=\frac{e^{2^n}}{n+1}\rightarrow\infty$.
Thus the radius of convergence of the Taylor series of $f$ at $t=0$
is zero. Now since $f$ is periodic with period $\pi$ we conclude
that the radius of convergence of the Taylor series of $f$ at
$t=k\pi$ is zero for all $k\in\Z$.

Next for any integer $N$ we set
$g_N(t)=\sum_{m=0}^{N}\frac{1}{m!}e^{i2^mt}$ and
$h_N(t)=\sum_{m=N+1}^{\infty}\frac{1}{m!}e^{i2^mt}$. We have
$f=g_N+h_N$. Clearly $g_N$ is everywhere analytic. $h_N$ is periodic
with period $\frac{\pi}{2^N}.$ Thus by a similar argument as above
we conclude that the radius of convergence of the Taylor series of
$h_N$ at $t=\frac{k\pi}{2^N}$ is zero for all $k\in\Z$. Thus the
radius of convergence of the Taylor series of $f$ at
$t=\frac{k\pi}{2^N}$ is zero for all $k\in\Z$, too. Since the set of
all numbers $\frac{k\pi}{2^N}, k\in \Z,N\in\N$ is dense in $\R$, we
conclude that the radius of convergence of the Taylor series of $f$
at any $t\in\R$ is zero. Thus by the theorem \ref{analytic2} part
(i), $\hat{f}(t)$ diverges for all $t\ne0$. \eproof
\begin{eg}
The function
\eqn{}{f(t):=\sum_{m=1}^{\infty}\frac{1}{m!}e^{i2^{-m}t}}is analytic
at $t=0$ whose convergence radius is infinity. Thus $\hat{f}(t)$
converges for all $t$ to $f(0)$.
\end{eg}
\proof The proof is similar to the previous example.\eproof

\textbf{Open Questions} 1. For the function
$$f(t):=e^{\frac{-1}{t^2}},~t\ne0,~f(0)=0,$$ does $\hat{f}(t)$
converge?

2. Is there non-analytic functions $f$ such that the series
$\hat{f}$ is point-wise or uniformly convergent and among such
functions if there is any, is there any function such that the sum
of the series is non-constant?

3. Verify the convergence of $\hat{f}(t)$ in the example \ref{u} for
$t\ne\frac{(2m+1)\ell}{2^n}$.

4. If we define a linear differential operator of infinite
order\eqn{}{f\mapsto
f(0)-\sum_{n=0}^{\infty}\frac{(-1)^n}{n!}t^nf^{(n)}(t)} then in
above we showed that analytic functions around origin are contained
in the space of eigenfunctions of the zero eigenvalue of this
operator. Now the question arises that: are there nonzero
eigenvalues for this operator?

\end{document}